\input amstex
\documentstyle{amsppt}
\topmatter
\title{Elementary Proofs of Identities for Schur Functions and Plane
Partitions}\endtitle
\author{David M.\ Bressoud}\endauthor
\date{August 24, 1998}\enddate
\dedicatory To George Andrews on the occasion of his 60th birthday \enddedicatory
\thanks Copyright to this work is retained by the
author. Permission is granted for the noncommercial
       reproduction of the complete work for educational or research purposes, and for the use of figures,
       tables and short quotes from this work in other books or journals, provided a full bibliographic
       citation is given to the original source of the material. \endthanks
\abstract  We use elementary methods to prove product formulas for sums of restricted
classes of Schur functions. These imply known identities for the generating function
for symmetric plane partitions with even column height and for the generating
function for symmetric plane partitions with an even number of angles at each level.
\endabstract
\address  Dept.\ of Mathematics \& Computer Science, Macalester College,
St.\ Paul, MN 55105, USA 
\endaddress
\email bressoud\@macalester.edu \endemail
\endtopmatter
\define\pf{\noindent {\bf Proof:} \ }
\define\inv{{\Cal{I}}}

\define\cS{{\Cal{S}}}
\define\calB{{\Cal{B}}}
\define\het{\text{Ht}}
\document

\head 1. Introduction \endhead

\noindent By a {\bf plane partition}, we mean a finite set, $\Cal{P}$, of lattice
points with positive integer coefficients, $\{(i,j,k)\} \subseteq {\Bbb N}^3$, with
the property that if $(r,s,t) \in \Cal{P}$ and $1 \leq i \leq r,\ 1 \leq j \leq s,\ 1
\leq k \leq t$, then $(i,j,k)$ must also be in $\Cal{P}$. A plane partition is {\bf
symmetric} if $(i,j,k) \in \Cal{P}$ if and only if $(j,i,k)\in \Cal{P}$. The {\bf
height of stack} $(i,j)$ is the largest value of $k$ for which there exists a point
$(i,j,k)$ in the plane partition. A plane partition is {\bf column strict} if the
height of stack $(i,j)$ is strictly less than the height of stack $(i-1,j)$
whenever $i \geq 2$ and $(i,j,1)$ is in the plane partition. 

Symmetric plane partitions were studied by P.\ A.\ MacMahon \cite{12} who
conjectured in 1898 that the generating function for symmetric plane partitions with
$1 \leq i, j \leq n $ and $1 \leq k \leq m$ is
$$
\prod_{i=1}^n \frac{1-q^{m+2i-1}}{1-q^{2i-1}} \prod_{1 \leq i < j \leq
n}\frac{1-q^{2(m+i+j-1)}}{1-q^{2(i+j-1)}}.
$$

This was proven independently by Andrews \cite{1} and Macdonald \cite{11}. As
shown by Andrews \cite{2}, this is equivalent to the Bender-Knuth conjecture
\cite{3}, that the generating function for column strict plane partitions with $1
\leq i,k \leq n$, $1 \leq k \leq m$ is
$$ \prod_{1 \leq i \leq j \leq n} {1-q^{m+i+j-1} \over 1-q^{i+j-1}}. $$

Both of these generating functions are consequences of the following theorem of
Macdonald \cite{11}, the first when we set
$x_i = q^{2n-2i+1}$ and the second when we set $x_i = q^{n+1-i}$.

\proclaim{Theorem I (Macdonald)} For positive integers $m$ and $n$,
$$
\sum_{\lambda\subseteq \{m^n\}} s_{\lambda}(x_1,\ldots,x_n) =
\frac{\det(x_i^{j-1}-x_i^{m+2n-j})}{\prod_{i=1}^n (1-x_i) \prod_{1 \leq i < j
\leq n} (x_ix_j-1)(x_i-x_j)}.\tag1.1
$$
The sum is over all partitions into at most $n$ parts, each of which is less than or
equal to $m$.
\endproclaim

I gave an elementary proof of Macdonald's identity in \cite{4}. D\'esarm\'enien
\cite{7} and Stembridge \cite{15} found a similar theorem where the sum on
the left is over partitions into even parts. D\'esarm\'enien \cite{8} has also found
the generalization in which any number of odd parts are specified.

\proclaim{Theorem II (D\'esarm\'enien-Stembridge)} For positive even integer $m$ and
positive integer $n$,
$$
\sum_{\lambda\subseteq \{m^n\} \atop \lambda\ \text{even}}
s_{\lambda}(x_1,\ldots,x_n) =
\frac{\det(x_i^{j-1}-x_i^{m+2n+1-j})}{\prod_{i=1}^n(1-x_i^2) \prod_{1 \leq i < j
\leq n} (x_ix_j-1)(x_i-x_j)}.\tag1.2
$$
\endproclaim

This theorem has two corollaries that were found by D\'esarm\'enien and Stembridge and,
independently, Proctor \cite{14}. The $q=1$ case was first
discovered by DeSainte-Catherine and Viennot \cite{6}. The generating function for
symmetric plane partitions with $1 \leq i,j \leq n$ and $1 \leq k \leq m$ where
$m$ is even and every stack has even height is given by
$$ \prod_{i=1}^n \frac{1-q^{m+2i}}{1-q^{2i}} \prod_{1 \leq i < j \leq n} \frac{1 -
q^{2(m+i+j)}}{1-q^{2(i+j)}}.
$$
The generating function for column strict plane partitions with $1 \leq i,k \leq n$,
$1 \leq j \leq m$ ($m$ even), and all rows of even length is
$$ \prod_{i=1}^n \frac{1-q^{m+2i}}{1-q^{2i}} \prod_{1 \leq i < j \leq n} \frac{1 -
q^{m+i+j}}{1-q^{i+j}}.
$$

Okada \cite{13} has proven the following companion using his minor summation
formula. Krattenthaler \cite{10} has used the special orthogonal
tableaux of Lakshmibai, Musili, and Seshadri to generalize this result to one in
which the number of columns of odd length is specified.

\proclaim{Theorem III (Okada)} For positive integer $m$ and positive even
integer
$n$,
$$
\sum_{\lambda\subseteq \{m^n\} \atop \lambda' \text{even}}
s_{\lambda}(x_1,\ldots,x_n) = {1\over 2}\,
\frac{\det(x_i^{j-1}-x_i^{m+2n-1-j})+\det(x_i^{j-1}+x_i^{m+2n-1-j})}{\prod_{1 \leq i
< j \leq n} (x_ix_j-1)(x_i-x_j)},\tag1.3
$$
where $\lambda'$ is the partition conjugate to $\lambda$. In other words the sum on
the left is over partitions with even column lengths.
\endproclaim

This has the following corollary when $x_i = q^{2n-2i+1}$.

\proclaim{Corollary} The generating function for symmetric plane partitions, $1 \leq
i, j \leq n$ ($n$ even) and $1 \leq k \leq m$, such that for each $k$ there are an
even number of lattice points of the form $(i,i,k)$ is given by
$$ {1 \over 2} \left(\prod_{i=0}^{n-1} (1-q^{m+2i}) + \prod_{i=0}^{n-1} (1+q^{m+2i})
\right) \prod_{1 \leq i < j \leq n} \frac{1-q^{2(m+i+j-2)}}{1-q^{2(i+j-1)}}.
$$
\endproclaim

The generating functions that are derived from Theorems~I, II, and III have
particularly nice formulations. We define $\calB(n,n,m) =
\{(i,j,k)\,|\,1 \leq i, j \leq n,\ 1 \leq k \leq m\}$ and $\calB(n,n,m)/\cS_2$
to be the set of orbits of $\calB(n,n,m)$ under transposition of the
first two coordinates. For $\eta \in \calB(n,n,m)/\cS_2$, we define $\het(\eta) =
i+j+k-2$ where $(i,j,k)$ is any one element of $\eta$. An {\bf orbit counting}
generating function is the sum in which each plane partition is weighted by $q$ to
the number of orbits. 

The generating function for symmetric plane partitions in $\calB(n,n,m)$ is given by
$$ \prod_{\eta\in\calB(n,n,m)/\cS_2}
\frac{1-q^{|\eta|(1+\het(\eta))}}{1-q^{|\eta|\,\het(\eta)}}.$$
The orbit counting generating function for symmetric plane partitions in
$\calB(n,n,m)$ is given by
$$ \prod_{\eta\in\calB(n,n,m)/\cS_2}
\frac{1-q^{1+\het(\eta})}{1-q^{\het(\eta)}}.$$
The generating function for symmetric plane partitions with even stack height in
$\calB(n,n,m)$ ($m$ even) is given by
$$ \prod_{\eta\in\calB(n,n,m)/\cS_2}
\frac{1-q^{|\eta|(2+\het(\eta))}}{1-q^{|\eta|(1+\het(\eta))}}.$$
The orbit counting generating function for symmetric plane partitions with even
stack height in $\calB(n,n,m)$ ($m$ even) is given by
$$ \prod_{\eta\in\calB(n,n,m)/\cS_2}
\frac{1-q^{2+\het(\eta})}{1-q^{1+\het(\eta)}}.$$
The generating function for symmetric plane partitions in $\calB(n,n,m)$ ($n$
even) such that for each $k$, $1 \leq k \leq m$, there are an even number of points
of the form $(i,i,k)$ is given by
$$ \prod_{\eta\in\calB(n,n,m-1)/\cS_2}
\frac{1-q^{|\eta|(1+\het(\eta))}}{1-q^{|\eta|\,\het(\eta)}} \sum_{S\subseteq
\{(i,i,m)\,|\,1 \leq i \leq n\} \atop |S| \text{ even}} \prod_{\eta\in
S}q^{\het(\eta)}.$$
There is a formula given by Krattenthaler (equation(7.15) in \cite{10}) for the
corresponding orbit counting generating function. It is not as readily stated in
terms of orbits.

In section~2, we shall warm up to the proof of Theorems~II and III with a general
result that includes the limiting cases of Theorems~I, II, and III. It was first
proved by Ishikawa and Wakayama \cite{9} using Okada's minor-summation formula of
Pfaffians.

\proclaim{Theorem IV (Ishikawa and Wakayama)} For any positive integer $n$, we have
that
$$ \sum_{\lambda} f_{\lambda}(t,v) s_{\lambda}(x_1,\ldots,x_n) = \prod_{i=1}^n
\frac{1}{(1-tx_i)(1-vx_i)} \prod_{1 \leq i < j \leq n} \frac{1}{1-x_ix_j},\tag1.4 $$
where we let $a_j$ be the number of columns of length $j$ in $\lambda$
(equivalently, the number of parts of size $j$ in $\lambda'$) and
$$ f_{\lambda}(t,v) = \prod_{j\ \text{odd}}
\frac{v^{a_j+1}-t^{a_j+1}}{v-t}\prod_{j\ \text{even}} \frac{1-(tv)^{a_j+1}}{1-tv} .
$$
\endproclaim

We note that
$$ \align f(0,1) &= 1, \\
f(1,-1) &= \left\{ \aligned 0 &\quad \text{if any $a_j$ is odd}, \\ 1 &\quad
\text{otherwise}, \endaligned \right. \\
f(0,0) &= \left\{ \aligned 0 &\quad \text{if any $a_j$ is positive for any odd $j$},
\\ 1 &\quad
\text{otherwise}. \endaligned \right.
\endalign $$
Theorem~IV implies the following Littlewood formulas
(\cite{11}, examples 4 and 5 in section I.5):
$$ \align
\sum_{\lambda} s_{\lambda}(x_1,\ldots,x_n) &= \prod_{i=1}^n
\frac{1}{1-x_i} \prod_{1 \leq i < j \leq n} \frac{1}{1-x_ix_j}, \tag1.5 \\
\sum_{\lambda \text{ even}} s_{\lambda}(x_1,\ldots,x_n) &=
\prod_{i=1}^n
\frac{1}{1-x_i^2} \prod_{1 \leq i < j \leq n} \frac{1}{1-x_ix_j}, \tag1.6 \\
\sum_{\lambda' \text{ even}} s_{\lambda}(x_1,\ldots,x_n) &=
\prod_{1 \leq i < j \leq n} \frac{1}{1-x_ix_j}. \tag1.7
\endalign $$

In section~3, we shall give the proof of Theorem~III as well as a new proof of
Theorem~II. Section~4 will show the derivation of the generating function for
symmetric plane partitions with an even number of lattice points of the form
$(i,i,k)$ for each $k$. With the exceptions of Lemmas 1 and 2, the results presented
in this paper are not new. The proofs, however, are considerably simpler than those
that have been given before.

\head 2. Proof of Theorem IV \endhead

\proclaim{Lemma 1} For any positive integer $n$ we have that
$$ \multline x_1\cdots x_n \sum_{k=1}^n x_k^{-1}(1-tx_k)(1-vx_k) \prod_{i=1 \atop i
\neq k}^n
\frac{1-x_ix_k}{x_i-x_k} \\
 = \left\{ \aligned &(1-tx_1\cdots x_n)(1-vx_1\cdots x_n),
\quad  \text{\rm if $n$ is odd}, \\ &(1-x_1\cdots x_n)(1-tvx_1\cdot x_n),
\quad  \text{\rm if $n$ is even}.\endaligned\right. \endmultline\tag2.1$$
\endproclaim

\pf This lemma is correct for $n=1$. We assume that it is correct with $n-1$
variables and proceed by induction. If we multiply both sides of equation~(2.1) by
$\prod_{1 \leq i < j \leq n}(x_i-x_j)$, each side is an alternating polynomial in
$x_1, \ldots, x_n$. It follows that both sides of equation~(2.1) are symmetric
polynomials that are quadratic in each of the variables $x_1$ through $x_n$. We
only need to show that they agree at three values of $x_1$. Both polynomials are 1
when $x_1=0$. When $x_1 = t^{-1}$, the polynomial on the left is equal to 
$$ \split & t^{-1}x_2\cdots x_n \sum_{k=2}^n x_k^{-1}(1-tx_k)(1-vx_k)\,
\frac{1-t^{-1}x_k}{t^{-1}-x_k} \prod_{i=2
\atop i
\neq k}^n \frac{1-x_ix_k}{x_i-x_k} \, \frac{1-t^{-1}x_k}{t^{-1}-x_k} \\
& \qquad = x_2\cdots x_n \sum_{k=2}^n x_k^{-1}(1-t^{-1}x_k)(1-vx_k) \prod_{i=2
\atop i
\neq k}^n \frac{1-x_ix_k}{x_i-x_k}\\
& \qquad = \left\{ \aligned &(1-x_2\cdots x_n)(1-t^{-1}vx_2\cdots x_n),
\quad  \text{if } n \text{ is odd}, \\ &(1-t^{-1}x_2\cdots x_n)(1-vx_2\cdot
x_n),
\quad  \text{if } n \text{ is even}.\endaligned\right. \endsplit \tag2.2$$
Similarly, the two polynomials agree at $x_1 = v^{-1}$. \quad $\square$

\proclaim{Lemma 2} For even positive integer $n$ we have that
$$ (x_1\cdots x_n)^2 \sum_{k=1}^n\sum_{l=1 \atop l \neq k}^n
x_k^{-2}x_l^{-1} \prod_{i=1 \atop i \neq k}^n \frac{1-x_ix_k}{x_i-x_k} \prod_{i=1
\atop i \neq k,l}^n \frac{1-x_ix_l}{x_i-x_l}
= 1-x_1\cdots x_n. \tag2.3 $$
\endproclaim

\pf This follows from lemma~1 with $t=v=0$, summing first over $l$ and then over
$k$. \quad $\square$

\medskip

\noindent{\bf Proof of Theorem~IV:} When $n=1$, the left side of equation~(1.4)
is
$$ \sum_{k=0}^{\infty} \frac{v^{k+1}-t^{k+1}}{v-t}\, x^k = \frac{1}{(1-vx)(1-tx)}. $$
We proceed by induction and assume that the equation is valid for $n-1$ variables.
We rewrite equation~(1.4) as
$$\multline \sum_{\lambda}\sum_{\sigma\in\cS_n}
(-1)^{\inv(\sigma)} f_{\lambda}(t,v) \prod_{i=1}^n
x_i^{\lambda_{\sigma(i)}+n-\sigma(i)}(1-tx_i)(1-vx_i)\prod_{1\leq i < j
\leq n}(1-x_ix_j) \\ = \prod_{1 \leq i < j \leq n} (x_i-x_j),
 \endmultline \tag2.3
$$
where $\inv(\sigma)$ is the inversion number of $\sigma$. We shall prove that the
left side is equal to the Vandermonde determinant. 

We take the double summation and first sum over the possible values of $\lambda_n$
and $k = \sigma^{-1}(n)$. We let $\tau$ be the restriction of $\sigma$ to
$\{1,\ldots n\}\backslash\{k\}$. If we subtract $\lambda_n$ from each of the parts
in
$\lambda$, we are left with a partition, $\mu$, into at most $n-1$ parts. We have
that
$f_{\lambda}(t,v) = c_{\lambda_n} f_{\mu}(t,v)$ where $c_{\lambda_n}$ is
$(v^{\lambda_n+1}-t^{\lambda_n+1})/(v-t)$ if $n$ is
odd, $(1-(vt)^{\lambda_n+1})/(1-vt)$ if $n$ is even. The left side of equation~(2.3)
is equal to
$$ \multline \sum_{\lambda_n=0}^{\infty} \sum_{k=1}^n
(-1)^{n-k}x_k^{-1}(1-tx_k)(1-vx_k)c_{\lambda_n} (x_1\cdots
x_n)^{\lambda_n+1}\prod_{i=1 \atop i \neq k}^n(1-x_ix_k) \\
\times \sum_{\mu,\tau} (-1)^{\inv(\tau)} f_{\mu}(t,v) \prod_{i=1 \atop i\neq k}^n
x_i^{\mu_{\tau(i)} + n-1 - \tau(i)}(1-tx_i)(1-vx_i) \prod_{1 \leq i < j \leq n
\atop i,j \neq k} (1-x_ix_j).
\endmultline $$

We use our induction hypothesis to rewrite this as
$$ \prod_{1 \leq i < j \leq n}(x_i-x_j) \sum_{\lambda_n=0}^{\infty} \sum_{k=1}^n
x_k^{-1}(1-tx_k)(1-vx_k)c_{\lambda_n} (x_1\cdots
x_n)^{\lambda_n+1}\prod_{i=1 \atop i \neq k}^n\frac{1-x_ix_k}{x_i-x_k}. $$
By Lemma~1, the double sum is equal to 1. \quad $\square$

\head 3. Proof of Theorems II and III \endhead

\noindent The proofs of Theorems~II and III are similar in structure to the proof
of Theorem~IV, just more complicated in detail.

\medskip

\noindent{\bf Proof of Theorem~II:} We verify that this theorem is correct for
$n=1$ and proceed by induction on the number of variables. We shall prove this
theorem in the form
$$ \multline \sum_{\lambda\subseteq \{m^n\} \atop \lambda \text{ even}}
\det(x_i^{\lambda_j+n-j}) \prod_{i=1}^n(1-x_i^2) \prod_{1 \leq i < j \leq
n}(x_ix_j-1) \\ = \sum_{\sigma\in\cS_n}\sum_{S \subseteq \{1,\ldots,n\}}
(-1)^{\inv(\sigma) + |S|} \prod_{i\in S}x_i^{m+2n+1-\sigma(i)} \prod_{i\notin
S}x_i^{\sigma(i)-1}, \endmultline
\tag3.1 $$
where $m$ is an even integer. 

As in the proof of Theorem~IV, we expand the left side as a sum over
partitions, $\lambda$, and permutations, $\sigma$. We then sum separately over
$\lambda_n$ which must now be even, $\lambda_n = 2t$, and over $k = \sigma^{-1}(n)$,
leaving $\mu$, the partition obtained from $\lambda$ when $\lambda_n$ is subtracted
from each part, and $\tau$, the restriction of $\sigma$ to
$\{1,\ldots,n\}\backslash\{k\}$. We then apply our induction hypothesis. The left
side of equation~(3.1) becomes
$$ \align & \sum_{t=0}^{m/2} \sum_{k=1}^n (-1)^{n-k} x_k^{-1}(1-x_k^2) (x_1 \cdots
x_n)^{2t+1} \prod_{i=1 \atop i \neq k}^n (x_ix_k-1) \\ & \quad \times \sum_{\mu,
\tau} (-1)^{\inv(\tau)} \prod_{i=1 \atop i \neq
k}^nx_i^{\mu_{\tau(i)}+n-1-\tau(i)}(1-x_i^2) \prod_{1 \leq i < j \leq n \atop i,j
\neq k}(x_ix_j-1) \\
& \quad = \sum_{t=0}^{m/2} \sum_{k=1}^n (-1)^{n-k} x_k^{-1}(1-x_k^2) (x_1 \cdots
x_n)^{2t+1} \prod_{i=1 \atop i \neq k}^n (x_ix_k-1) \\
& \qquad \times \sum_{\sigma\in \cS_{n-1}} \sum_{S \subseteq
\{1,\ldots,n\}\backslash\{k\}} (-1)^{\inv(\sigma) + |S|} \prod_{i\in S}
x_i^{m-2t+2(n-1)+1 - \sigma(i)} \prod_{i\in\overline{S}}x_i^{\sigma(i)-1},
\endalign$$
where $\cS_{n-1}$ is the set of 1--1 mappings from $\{1,\ldots,n\}\backslash\{k\}$
to $\{1,\ldots,n-1\}$ and $\overline{S}$ is the complement of $S$ in
$\{1,\ldots,n\}\backslash\{k\}$.

We simplify this summation and then sum over $\sigma\in \cS_{n-1}$ and over $t$.
The left side of equation~(3.1) is equal to
$$ \align
& \sum_{t=0}^{m/2} \sum_{k=1}^n  \sum_{\sigma\in
\cS_{n-1}} \sum_{S \subseteq
\{1,\ldots,n\}\backslash\{k\}} (-1)^{n-k+\inv(\sigma) + |S|} x_k^{-1}(1-x_k^2)
\prod_{i\notin S} x_i^{2t+1}
\prod_{i=1 \atop i \neq k}^n (x_ix_k-1) \\
& \quad \times \prod_{i\in S}x_i^{m+2n-\sigma(i)}
\prod_{i\in\overline{S}}x_i^{\sigma(i)-1} \\
& \quad = \sum_{k=1}^n\sum_{S \subseteq \{1,\ldots,n\}\backslash\{k\}}
(-1)^{n-k+|S|} x_k^{-1}(1-x_k^2) \frac{1-\prod_{i\notin
S}x_i^{m+2}}{1-\prod_{i\notin S}x_i^2} \prod_{i=1 \atop i \neq k}^n (x_ix_k-1) \\
& \qquad \times \prod_{i\notin S}x_i \prod_{i\in S}x_i^{m+2n-1} (-1)^{\left(n-1\atop
2\right)}
\prod_{1
\leq i < j
\leq n \atop i,j \neq k} (x_i^{\epsilon_i}-x_j^{\epsilon_j}),
\endalign $$
where $\epsilon_i$ is $-1$ if $i \in S$ and $+1$ if $i \notin S$.

We reverse the order of summation so that we first sum over all proper subsets of
$\{1,\ldots,n\}$ and then over all $k \notin S$. For each $i\in S$, we rewrite
$x_ix_k-1$ as $-x_i(x_i^{\epsilon_i}-x_k^{\epsilon_k})$ if $i < k$, and rewrite it
as  $x_i(x_k^{\epsilon_k}-x_i^{\epsilon_i})$ if $i>k$. The left side of
equation~(3.1) has become
$$ \align
&(-1)^{\left(n\atop 2\right)}  \sum_{S\subset\{1,\ldots,n\}}
(-1)^{|S|} \prod_{i\in S}x_i^{m+2n} \frac{1-\prod_{i\notin
S}x_i^{m+2}}{1-\prod_{i\notin S}x_i^2} \prod_{1 \leq i < j \leq
n}(x_i^{\epsilon_i}-x_j^{\epsilon_j}) \\
& \quad \times \prod_{i\notin S}x_i \sum_{k\notin S}x_k^{-1}(1-x_k^2)
\prod_{i\notin S \atop i \neq k} \frac{1-x_ix_k}{x_i-x_k}.
\endalign $$

By Lemma~1, the second line is equal to $1-\prod_{i\notin S} x_i^2$ which cancels
with the factor in the denominator. We now expand the Vandermonde product. The left
side of equation~(3.1) is equal to
$$ \align &
\sum_{S\subset\{1,\ldots,n\}}\sum_{\sigma\in \cS_n} (-1)^{\inv(\sigma) + |S|}
\prod_{i\in S}x_i^{m+2n+1-\sigma(i)} \prod_{i\notin S} x_i^{\sigma(i)-1} \\
& \quad - \sum_{S\subset\{1,\ldots,n\}}\sum_{\sigma\in \cS_n} (-1)^{\inv(\sigma) + |S|}
\prod_{i\in S}x_i^{m+2n+1-\sigma(i)} \prod_{i\notin S} x_i^{m+\sigma(i)+1}.
\endalign $$

We use the fact that
$$ \multline
 \sum_{S\subseteq\{1,\ldots,n\}}\sum_{\sigma\in \cS_n} (-1)^{\inv(\sigma) + |S|}
\prod_{i\in S}x_i^{m+2n+1-\sigma(i)} \prod_{i\notin S} x_i^{m+\sigma(i)+1} \\
 = \det\left(x_i^{m+j+1} - x_i^{m+2n+1-j}\right) = 0,
\endmultline $$
to replace
$$ - \sum_{S\subset\{1,\ldots,n\}}\sum_{\sigma\in \cS_n} (-1)^{\inv(\sigma) + |S|}
\prod_{i\in S}x_i^{m+2n+1-\sigma(i)} \prod_{i\notin S} x_i^{m+\sigma(i)+1} $$
by 
$$\sum_{\sigma\in\cS_n}(-1)^{\inv(\sigma)+n} \prod_{i=1}^n x_i^{m+2n+1-\sigma(i)}. $$
This puts the left side of equation~(3.1) in the desired form. \quad $\square$

\medskip

\noindent{\bf Proof of Theorem~III:} We begin by rewriting the identity to be
proven as
$$ \multline
 \sum_{\lambda\subseteq\{m^n\}\atop \lambda' \text{ even}}
\sum_{\sigma\in\cS_n}(-1)^{\inv(\sigma)}\prod_{i=1}^n
x_i^{\lambda_{\sigma(i)}+n-\sigma(i)} \prod_{1 \leq i < j \leq n} (x_ix_j-1) \\
= \sum_{\sigma\in\cS_n} \sum_{S\subseteq \{1,\ldots,n\} \atop |S| \text{ even}}
(-1)^{\inv(\sigma)} \prod_{i\in S} x_i^{m+2n-1-\sigma(i)}\prod_{i\notin
S}x_i^{\sigma(i)-1}.
\endmultline  \tag3.2 $$

We again proceed by induction. For this theorem, we need to identify both
$\sigma^{-1}(n)$ and $\sigma^{-1}(n-1)$. We form $\mu$ by subtracting $\lambda_n$
from each part. Since each column has even length, $\mu$ has at most $n-2$ parts.
The left side of equation~(3.2) is equal to
$$\align
& \sum_{\lambda_n=0}^m \sum_{1 \leq k < l \leq
n}(-1)^{n-k+n-l} (x_k^{-2}x_l^{-1}-x_k^{-1}x_l^{-2}) (x_1\cdots
x_n)^{\lambda_n+2} \\
& \quad \times (x_kx_l-1)\prod_{i=1 \atop i \neq k,l}^n(x_ix_k-1)(x_ix_l-1) \\
& \quad \times \sum_{\mu, \tau}(-1)^{\inv(\tau)} \prod_{i=1 \atop i \neq k,l}^n
x_i^{\mu_{\tau(i)}+n-2-\tau(i)} \prod_{1 \leq i < j \leq n \atop i,j \neq k,l}
(x_ix_j-1).
\endalign$$

We apply our induction hypothesis to the inner sum and then sum over $\lambda_n$ and
$\sigma\in\cS_{n-2}$, the set of 1--1 mappings from
$\{1,\ldots,n\}\backslash\{k,l\}$ to $\{1,\ldots,n-2\}$. The left side of
equation~(3.2) becomes
$$ \align
& \sum_{1 \leq k < l \leq n}\, \sum_{S\subseteq\{1,\ldots,n\}\backslash\{k,l\}\atop
|S| \text{ even}} (-1)^{k+l}(x_k^{-2}x_l^{-1}-x_k^{-1}x_l^{-2}) 
\frac{1-\prod_{i\notin S}x_i^{m+1}}{1-\prod_{i\notin S}x_i} \\
& \quad \times (x_kx_l-1)\prod_{i=1 \atop i \neq k,l}^n(x_ix_k-1)(x_ix_l-1) \\
& \quad \times \prod_{i\notin S}x_i^2 \prod_{i\in S} x_i^{m+2n-4} (-1)^{\left(n-2
\atop 2\right)} \prod_{1 \leq i < j \leq n \atop i,j \neq k,l}(x_i^{\epsilon_i} -
x_j^{\epsilon_j}).
\endalign $$
Again we have $\epsilon_i = -1$ if $i\in S$ and $+1$ if $i \notin S$.

For $i\in S$, we rewrite $(x_ix_k-1)(x_ix_l-1)$ as
$x_i^2(x_i^{\epsilon_i}-x_k^{\epsilon_k})(x_i^{\epsilon_i}-x_l^{\epsilon_l})$. We
then interchange the sum on $S$, proper subsets of $\{1,\ldots,n\}$ with even
cardinality, and the sum on $k$ and $l$ which now must lie in the complement of
$S$. The left side of equation~(3.2) is equal to
$$ \align
& (-1)^{\left(n \atop 2\right)} \sum_{S \subset\{1,\ldots,n\}\atop |S| \text{ even}}
\prod_{i\in S}x_i^{m+2n-2} \frac{1-\prod_{i\notin S}x_i^{m+1}}{1-\prod_{i\notin
S}x_i} \prod_{1 \leq i < j \leq n}(x_i^{\epsilon_i}-x_j^{\epsilon_j}) \\
& \quad \times \prod_{i\notin S}x_i^2 \sum_{1 \leq k < l \leq n \atop k,l \notin
S} (x_k^{-2}x_l^{-1}-x_k^{-1}x_l^{-2}) \frac{x_kx_l-1}{x_k-x_l} \prod_{i\notin S
\atop i \neq k,l} \frac{(1-x_ix_k)(1-x_ix_k)}{(x_i-x_k)(x_i-x_l)}.
\endalign $$
The second line of this expression is equal to
$$ \prod_{i\notin S}x_i^2\sum_{k\notin S} \sum_{l\notin S \atop l \neq
k}x_k^{-2}x_l^{-1} \prod_{i \notin S \atop i \neq k} \frac{1-x_ix_k}{x_i-x_k}
\prod_{i\notin S \atop i \neq k,l} \frac{1-x_ix_l}{x_i-x_l}.
$$
By Lemma~2, this is equal to $1-\prod_{i\notin S}x_i$ which cancels with the factor
in the denominator.

As in the proof of Theorem~II, we replace the Vandermonde product by the sum over
permutations. The left side of equation~(3.2) becomes
$$ \align
& \sum_{S\subset\{1,\ldots,n\}\atop |S| \text{ even}} \sum_{\sigma\in\cS_n}
(-1)^{\inv(\sigma)} \prod_{i\in S} x_i^{m+2n-\sigma(i)-1} \prod_{i\notin
S}x_i^{\sigma(i)-1} \\
& \quad - \sum_{S\subset\{1,\ldots,n\}\atop |S| \text{ even}} \sum_{\sigma\in\cS_n}
(-1)^{\inv(\sigma)} \prod_{i\in S} x_i^{m+2n-\sigma(i)-1} \prod_{i\notin
S}x_i^{m+\sigma(i)}.
\endalign $$
We now observe that
$$ 
\sum_{S\subseteq\{1,\ldots,n\}\atop |S| \text{ even}} \sum_{\sigma\in\cS_n}
(-1)^{\inv(\sigma)} \prod_{i\in S} x_i^{m+2n-\sigma(i)-1} \prod_{i\notin
S}x_i^{m+\sigma(i)} =  0. $$
This is true because if we interchange the inverse images of $n$ and $n-1$ and change
whether or not each inverse image is in $S$, then we change the sign of the
inversion number but do not change the monomial. As a result, we have that
$$\multline
- \sum_{S\subset\{1,\ldots,n\}\atop |S| \text{ even}} \sum_{\sigma\in\cS_n}
(-1)^{\inv(\sigma)} \prod_{i\in S} x_i^{m+2n-\sigma(i)-1} \prod_{i\notin
S}x_i^{m+\sigma(i)} \\
= \sum_{\sigma\in\cS_n} (-1)^{\inv(\sigma)} \prod_{i=1}^n x_i^{m+2n-\sigma(i)-1}
\endmultline $$
The left side of equation~(3.2) is equal to the desired sum. \quad $\square$ 

\head 4. Consequence and Question \endhead

\noindent If we set $x_i = q^{2n-2i+1}$ in Theorem~III, the left side of
equation~(1.3) becomes the generating function for symmetric plane partitions with
$1 \leq i,j \leq n$, $1 \leq k \leq m$, such that for each $k$ there are an even
number of lattice points of the form $(i,i,k)$. The right side of equation~(1.3)
becomes
$$ \align & \prod_{i=1}^n q^{(2n-2i+1)(m+2n-2)/2} \sum_{\sigma\in \cS_n}
\sum_{S\subseteq\{1,\ldots,n\}\atop |S| \text{ even}} (-1)^{\inv(\sigma)}
\prod_{i=1}^n q^{\epsilon_i(2n+m-2\sigma(i))(2i-2n-1)/2} \\
& \quad \times \prod_{1 \leq i < j \leq n}
(q^{2n-2i+1}-q^{2n-2j+1})^{-1}(q^{4n-2i-2j+1}-1)^{-1},
\endalign
$$
where $\epsilon_i = -1$ if $i\in S$ and $+1$ if $i \notin S$. 

We combine the $B_n$ form of the Weyl denominator formula,
$$ \multline \sum_{\sigma\in \cS_n}
\sum_{S\subseteq\{1,\ldots,n\}} (-1)^{\inv(\sigma)+|S|}
\prod_{i=1}^n x_{\sigma(i)}^{\epsilon_i(2i-2n-1)/2} \\ = \prod_{i=1}^n x_i^{(1-2n)/2}
(1-x_i) \prod_{1 \leq i < j \leq n} (x_i-x_j)(x_ix_j-1), \endmultline
$$
and the identity obtained when each $x_i$ is replace by $-x_i$:
$$ \multline \sum_{\sigma\in \cS_n}
\sum_{S\subseteq\{1,\ldots,n\}} (-1)^{\inv(\sigma)}
\prod_{i=1}^n x_{\sigma(i)}^{\epsilon_i(2i-2n-1)/2} \\ = \prod_{i=1}^n x_i^{(1-2n)/2}
(1+x_i) \prod_{1 \leq i < j \leq n} (x_i-x_j)(x_ix_j-1), \endmultline
$$
to derive the result that we need:
$$ \multline \sum_{\sigma\in \cS_n}
\sum_{S\subseteq\{1,\ldots,n\}\atop |S| \text{ even}} (-1)^{\inv(\sigma)}
\prod_{i=1}^n x_{\sigma(i)}^{\epsilon_i(2i-2n-1)/2} \\ = {1\over 2}\prod_{i=1}^n
x_i^{(1-2n)/2}
\left( \prod_{i=1}^n(1-x_i)+\prod_{i=1}^n(1+x_i)\right) \prod_{1 \leq i < j \leq n}
(x_i-x_j)(x_ix_j-1). \endmultline \tag4.1
$$
The corollary now follows directly.

It would be of interest to find the analogous formula for
$$ \sum_{\lambda \subseteq \{m^n\}} f_{\lambda}(t,v) s_{\lambda}(x_1,\ldots,x_n),
$$
although the form of it will certainly be much more complicated than anything
presented here.

\enddocument